\documentclass[12pt,fleqn]{article}
\usepackage{graphicx}

\usepackage{latexsym}

\usepackage{amsmath}
\usepackage{amsthm}
\usepackage{amssymb}
\usepackage{amsfonts}

\numberwithin{equation}{section}

\begin{document}

\newcommand{\rf}[1]{(\ref{#1})}
\newcommand{\rff}[2]{(\ref{#1}\ref{#2})}

\newcommand{\ba}{\begin{array}}
\newcommand{\ea}{\end{array}}

\newcommand{\be}{\begin{equation}}
\newcommand{\ee}{\end{equation}}

\newcommand{\const}{{\rm const}}
\newcommand{\ep}{\varepsilon}
\newcommand{\Cl}{{\cal C}}
\newcommand{\rr}{{\vec r}}
\newcommand{\ph}{\varphi}
\newcommand{\R}{{\mathbb R}}  
\newcommand{\N}{{\mathbb N}}
\newcommand{\Z}{{\mathbb Z}}

\newcommand{\e}{{\bf e}}

\newcommand{\m}{\left( \ba{r}}
\newcommand{\ema}{\ea \right)}
\newcommand{\mm}{\left( \ba{cc}}
\newcommand{\miv}{\left( \ba{cccc}}

\newcommand{\scal}[2]{\mbox{$\langle #1 \! \mid #2 \rangle $}}
\newcommand{\ods}{\par \vspace{0.3cm} \par}
\newcommand{\dis}{\displaystyle }
\newcommand{\sss}{\scriptscriptstyle}

\newcommand{\mc}{\multicolumn}
\newcommand{\no}{\par \noindent}

\newcommand{\sinc}{ {\rm sinc\,} }
\newcommand{\tanhc}{{\rm tanhc} }
\newcommand{\tanc}{{\rm tanc} }
\newcommand{\sech}{ {\rm sech\,} }

\newtheorem{prop}{Proposition}[section]
\newtheorem{Th}[prop]{Theorem}
\newtheorem{lem}[prop]{Lemma}
\newtheorem{rem}[prop]{Remark}
\newtheorem{cor}[prop]{Corollary}
\newtheorem{Def}[prop]{Definition}
\newtheorem{open}{Open problem}
\newtheorem{ex}[prop]{Example}
\newtheorem{exer}[prop]{Exercise}

\newenvironment{Proof}{\par \vspace{2ex} \par
\noindent \small {\it Proof:}}{\hfill $\Box$ 
\vspace{2ex} \par }

\title{\bf 
Improving the accuracy \\ of the AVF method }
\author{
 {\bf Jan L.\ Cie\'sli\'nski}\thanks{\footnotesize
 e-mail: \tt janek\,@\,alpha.uwb.edu.pl}
\\ {\footnotesize Uniwersytet w Bia{\l}ymstoku,
Wydzia{\l} Fizyki}
\\ {\footnotesize ul.\ Lipowa 41, 15-424
Bia{\l}ystok, Poland}
}

\date{}

\maketitle

\begin{abstract}

The Average Vector Field (AVF) method is a B-series scheme of the second order. As a discrete gradient method it preserves exactly the energy integral for any canonical Hamiltonian system. 
We present and discuss two locally exact and energy-preserving modifications of the AVF method:  AVF-LEX (of the third order) and AVF-SLEX (of the fourth order).   Applications to spherically symmetric potentials are given, including a compact explicit expression for the AVF scheme for the Coulomb-Kepler problem.

\end{abstract}

\no {\it PACS Numbers:} 45.10.-b; 02.60.Cb; 02.70.-c; 02.70.Bf 



\section{Introduction}

Discrete gradient methods were developed many years ago for molecular dynamics simulations \cite{LaG}. 
They preserve exactly the energy integral which is of  considerable advantage \cite{HLW}.  
More recently discrete gradient methods have been essentially developed by 
McLachlan, Quispel and their collaborators \cite{MQR2,MQTse,QM}. In particular, discrete gradient schemes  preserving any  first integrals for any ordinary differential equations were constructed \cite{QT}. 

A discrete gradient is a function of two vector variables ($\pmb{y}_{n}$,  $\pmb{y}_{n+1}$)  and has to satisfy (see \cite{MQR2}):
\be  \ba{l}
   \bar\nabla H (\pmb{y}_n,\pmb{y}_{n+1}) \cdot (\pmb{y}_{n+1} - \pmb{y}_n) = H (\pmb{y}_{n+1}) - H (\pmb{y}_n) \\[1ex]
\bar\nabla H (\pmb{y}_n,\pmb{y}_{n+1}) \rightarrow \nabla H (\pmb{y}) \quad \text{for} \quad  \pmb{y}_{n}  \rightarrow \pmb{y} \ , 
\ea \ee
where $H$ is a function of $\pmb{y}$ (we will confine ourselves to a particular  case when $H$ is a Hamiltonian).  We will consider only symmetric discrete gradients, i.e., such that  \  $ \bar\nabla H (\pmb{y}_n,\pmb{y}_{n+1}) =  \bar\nabla H (\pmb{y}_{n+1},\pmb{y}_{n}) $.  Discrete gradients are highly non-unique.  An important special case of the Average Vector Field  discrete gradient (AVF) was introduced in \cite{MQR2}. Later, the AVF scheme was identified as a B-series method \cite{QM}  which prompted intensive studies on energy-preserving B-series methods  \cite{CMMOQW,CMOQ,Ha}.  

 In this paper we derive simple explicit  formulas for the AVF discrete gradient in two important cases: the Coulomb-Kepler problem and a three-dimensional anharmonic oscillator, see Section~\ref{sec-AVF}.  We also  
propose two ``locally exact'' modifications  of the  AVF method: AVF-LEX and AVF-SLEX.  They preserve the energy integral and have much higher accuracy than the original AVF scheme. Although such modifications can be constructed for any Hamitonian and any discrete gradient (see \cite{Ci-locex-PLA}), in this paper we focus on discrete gradient schemes for   Hamiltonians of the form $H = \frac{1}{2} \pmb{p}^2  + V (\pmb{x})$, see Section~\ref{sec-Tay}. More general results can be found in Section~\ref{sec-B}. 

A numerical scheme for an ordinary differential equation  $\dot {\pmb y} = F ({\pmb y})$  is {\it locally exact} if there exists a sequence $(\pmb{\bar y}_n)$ such that  the linearization of the scheme 
around $\pmb{\bar y}_n$ is identical with the exact discretization of the differential equation linearized around $\pmb{\bar y}_n$, see  \cite{Ci-locex,CR-BIT}.   
In this paper we consider two types of locally exact modifications:   \ $\pmb{\bar y}_n = \pmb{y}_n$ \ (LEX) and  $\pmb{\bar y}_n = \frac{1}{2} \left( \pmb{y}_n + \pmb{y}_{n+1} \right)$  \ (SLEX). 
Our approach consists in considering  a class of non-standard modifications  of  a numerical scheme (compare \cite{Mic})  parameterized by some functions (e.g., matrices denoted by $\pmb{\delta}_n$). Requiring local exactness we can determine these functions (another possibility   is  proposed in \cite{CR-BIT,CR-CPC}).    
We point out that any linear ordinary differential equation with constant coefficients admits the exact discretization in an explicit form, see  \cite{Ci-oscyl,CR-ade,Mic,Potts}. 
The exact discretization of  linearized equations is known as the exponential Euler method \cite{Pope}. A similar notion (preservation of the linearization at all fixed points) appeared in \cite{MQTse}, see also \cite{CR-long}.

Locally exact modifications of the discrete gradient method for Hamiltonian systems  were studied  in \cite{Ci-locex,Ci-CMA,Ci-locex-PLA,CR-PRE,CR-BIT}.  Discrete gradient schemes were modified without changing their most important  property: the exact conservation of the energy integral. Locally exact modifications  improve the accuracy (especially in the neighbourhood of a stable equilibrium) although their order is not always higher. In the case of one degre of freedom the  symmetric discrete gradient scheme (of second order) has locally exact modifications of  third (LEX) and fourth (SLEX) order, respectively, \cite{CR-PRE,CR-BIT}. In the multidimensional case the order is usually unchanged, see \cite{Ci-CMA}. 
The AVF method turns out to be an exception: locally exact modifications, AVF-LEX and AVF-SLEX, are of third and fourth order, respectively. Numerical experiments confirm advantages of proposed modifications, see Section~\ref{sec-num}. Locally exact modifications increase the accuracy by several orders of magnitude.

\section{Average Vector Field method}
\label{sec-AVF}

We consider a Hamiltonian system of the form
\be    \label{F(y)}
    \pmb{\dot y} =  F (\pmb{y}) \ ,  \qquad    F (\pmb{y}) =  S \nabla H \ , \qquad     S = \mm 0 & I \\ -I &  0 \ema 
\ee
where $\pmb{y} \in \R^{2m}$, $H=H(\pmb{y})$ is a Hamiltonian,  and $I$ is  $m\times m$ unit matrix. The Average Vector Field method is defined by
\be  \label{AVF-y}
   \frac{\pmb{y}_{n+1} - \pmb{y}_n}{h_n} =  \int_0^1   F ( \pmb{y}_{n} +  \xi \Delta \pmb{y}_n ) d \xi  \  ,
\ee
where $h_n$ is a variable time step and $\Delta \pmb{y}_n := \pmb{y}_{n+1} - \pmb{y}_n$. 
The AVF method exactly preserves the energy integral  $H$, is a symmetric  B-series method of the second order, is affine-covariant and self-adjoint \cite{CMMOQW}. 

We denote $\pmb{y} = (\pmb{x}, \pmb{p}) = (x^1,\ldots,x^m,p^1,\ldots,p^m)$.  In this paper we often confine ourselves to the case   
\be  \label{H-p2}
    H = \frac{1}{2} \pmb{p}^2  + V (\pmb{x}) .
\ee 
Then
\be
  F (\pmb{y}) = \left(  \ba{c}  \pmb{p} \\ - V_{\pmb{x}} \ea  \right)  \equiv  \left( p^1, \ldots, p^m, - \frac{\partial V}{\partial x^1}, \ldots, - \frac{\partial V}{\partial x^m} \right)^T
\ee
and the AVF method  \rf{AVF-y} can be rewritten as 
\be  
\frac{\pmb{x}_{n+1} - \pmb{x}_n}{h_n} = \frac{ \pmb{p}_n + \pmb{p}_{n+1} }{2} , \qquad  
\frac{\pmb{p}_{n+1} - \pmb{p}_n}{h_n} = \int_0^1  \pmb{f} ( \pmb{x}_n + \xi \Delta \pmb{x}_n ) d \xi ,
 \ee
where $\pmb{f} (\pmb{x}) := - V_{\pmb{x}} (\pmb{x})$ and $\Delta \pmb{x}_n := \pmb{x}_{n+1} - \pmb{x}_n$. Therefore, it is natural to define the AVF discrete gradient as  
\be
  \bar \nabla^{\sss  AVF} V \equiv  \bar \nabla^{\sss  AVF} V (\pmb{x}_n, \pmb{x}_{n+1}) := - \int_0^1  \pmb{f} ( \pmb{x}_n + \xi \Delta \pmb{x}_n ) d \xi \ .
\ee

In some special cases the integral defining the AVF discrete gradient  can be explicitly calculated. We will consider the important case of a spherically symmetric anharmonic oscillator and the Coulomb-Kepler potential.

\begin{lem}   \label{lem-Simpson}
The AVF discrete gradient for an anharmonic oscillator potential \ $V (r) =  \alpha  r^2 - \beta r^4 $ \ 
(where $r = |\pmb{x}|$) \ can be explicitly computed by the Simpson rule  
\be  \label{Simpson} 
     \bar \nabla^{\sss  AVF}  V =  - \frac{1}{6} \left( \pmb{f}  (\pmb{x}_n)  + 4 \pmb{f} \left(  {\textstyle \frac{\pmb{x}_{n}+\pmb{x}_{n+1}}{2} }  \right)  +  \pmb{f} (\pmb{x}_{n+1})      \right) 
\ee
where $\pmb{f} (\pmb{x}) = - V_{\pmb{x}} (\pmb{x})$. 
\end{lem}

\begin{Proof}  If $V(r) = \frac{1}{n} r^n$, then   $\pmb{f}(\pmb{x}) =-  r^{n-2} \pmb{x}$. Therefore
\be
     \bar \nabla^{\sss  AVF}   \frac{ r^n}{n}  =  \int_0^1  \left(  (1-t) \pmb{x}_n + t \pmb{x}_{n+1} )^2  \right)^{\frac{n-2}{2} }             ( (1-t) \pmb{x}_n + t  \pmb{x}_{n+1}  ) d t 
\ee
which, in principle, can be integrated in an elementary way for any $n$. In particular, 
\be
     \bar \nabla^{\sss  AVF}   \frac{ r^2}{2}  =  \int_0^1       ( (1-t) \pmb{x}_n + t \pmb{x}_{n+1} ) d t  = \frac{1}{2} \left(  \pmb{x}_n + \pmb{x}_{n+1}  \right) \ ,
\ee
and 
\be  \ba{l}   \dis
     \bar \nabla^{\sss  AVF}   \frac{ r^4}{4}  =  \int_0^1  \left(  (1-t) \pmb{x}_n + t \pmb{x}_{n+1} )^2  \right)             ( (1-t) \pmb{x}_n + t  \pmb{x}_{n+1}  ) d t   \\[2em] \dis
\quad   =  \left( r_n^2  + \frac{1}{3} r_{n+1}^2  + \frac{2}{3} \pmb{x}_n \cdot \pmb{x}_{n+1}  \right)  \pmb{x}_n  +  \left( r_{n+1}^2  + \frac{1}{3} r_n^2  + \frac{2}{3} \pmb{x}_n \cdot \pmb{x}_{n+1}  \right)  \pmb{x}_{n+1} \ .
\ea \ee
Now the equality \rf{Simpson} can be verified as an easy exercise. 
\end{Proof}

 In a particular case, when $f (\pmb{x}_{n+1})$ is parallel to $\pmb{x}_{n+1}$, Lemma~\ref{lem-Simpson} is a consequence of  the well known fact that the Simpson rule is exact for cubic integrands. 

General results of this kind, concerning polynomial Hamiltonians, are reported in \cite{CMMOQW}.  The non-polynomial cases are more difficult. An important example ($n=-1$) is treated below.

\begin{lem}
The AVF discrete gradient for the Coulomb-Kepler  potential  \ $V (r) = -  \frac{\kappa}{r} $ \ 
can be explicitly computed as 
\be  \label{Kep}
     \bar \nabla^{\sss AVF}  V  =  \frac{ \kappa (\pmb{\hat x}_n + \pmb{\hat x}_{n+1})}{r_n r_{n+1} + \pmb{x}_n \cdot \pmb{x}_{n+1}} 
\ee
\end{lem}

\begin{Proof}  We denote   $r_n = |\pmb{x}_n|$ and  
\be   \label{defy1}
\pmb{\xi_n} = \frac{ \pmb{x}_{n+1} - \pmb{x}_n}{r_n}  \ ,   \quad   \pmb{\hat x}_n = \frac{\pmb{x}_n}{r_n} \ , \quad   \pmb{\hat \xi}_n = \frac{\pmb{\xi}_n}{|\pmb{\xi}_n|} \ ,   \quad  \cos\alpha_n = \pmb{\hat x}_n \cdot \pmb{\hat \xi}_n  \ .
\ee
Then
\be  \ba{l} \dis  \label{calka}
     \bar \nabla^{\sss AVF} \left( -  \frac{ 1}{r} \right)  =  \frac{1}{r_{n}^2} \int_0^1 \frac{ (\pmb{\hat x}_n + t \pmb{\xi}_n) d t  }{ ((\pmb{\hat x}_n + t \pmb{\xi}_n)^2    )^{3/2} } =   \frac{1}{r_{n}^2} \int_0^1 \frac{ (\pmb{\hat x}_n + t \pmb{\xi}_n) d t  }{ ((t \xi_n + \cos\alpha_n  )^2 + \sin^2\alpha_n     )^{3/2} }  \\[2em] 
\dis
\qquad \qquad = \frac{1}{r_{n}^2 \sin^2 \alpha_n \xi_n } \int_{u_0}^{u_1} \frac{ (\pmb{\hat x}_n + (\sin\alpha_n \sinh u - \cos\alpha_n) \pmb{\hat \xi}_n)  du   }{ \cosh^2 u }
\ea \ee
where  $\xi_n = |\pmb{\xi}_n|$ and 
\be
   \sinh u = \frac{t \xi_n + \cos\alpha_n}{\sin\alpha_n} \ , \quad \cosh u\,   du = \frac{\xi_n d t}{\sin\alpha_n} \ .
\ee
Furthermore,  $u_0$ and $u_1$ are defined by
\be  \label{sh0}
  \sinh u_0 = \cot\alpha_n \ , \quad \cosh u_0 = \frac{1}{\sin\alpha_n} \ ,  \quad  \sinh u_1 = \frac{\xi_n + \cos\alpha_n}{\sin\alpha_n} \ . 
\ee
It is convenient to consider the triangle formed by $\pmb{x}_n, \pmb{x}_{n+1},  r_n \pmb{\xi}_n$. We define
\be  \label{altyl} 
  \cos\tilde\alpha_n = \pmb{\hat x}_{n+1} \cdot \pmb{\hat \xi}_n \ .
\ee
Then also
\be   \label{differ}
 \pmb{x}_{n+1} \cdot \pmb{x}_n = r_{n+1} r_n \cos(\tilde\alpha_n -\alpha_n) \ , 
\ee
and, by the sine rule, 
\be  \label{sinerule}
  r_{n+1} \sin\tilde\alpha_n = r_n \sin\alpha_n \ , \qquad   \xi_n \sin\tilde\alpha = \sin(\alpha_n-\tilde\alpha_n) \ . 
\ee
Therefore, using \rf{defy1}, \rf{altyl} and \rf{sinerule}, we get
\be  \label{sh1}
 \sinh u_1 = \frac{\pmb{\xi}_n \cdot \pmb{\hat \xi}_n + \cos\alpha_n}{\sin\alpha_n} = \frac{r_{n+1} \cos\tilde\alpha }{r_n \sin\alpha} = \cot\tilde\alpha \ .
\ee
Therefore, \rf{calka} yields
\be
 \bar \nabla^{\sss AVF}  \frac{ 1}{r}  = \frac{(\tanh u_0 - \tanh u_1) (\pmb{\hat x}_n - \cos\alpha_n \pmb{\hat \xi}_n ) + (\sech u_1 - \sech u_0) \sin\alpha_n \pmb{\hat \xi}_n }{r_n^2\sin^2\alpha_n \xi_n} 
\ee
and, using \rf{sh0} and \rf{sh1}, we obtain
\be
 \bar \nabla^{\sss AVF} \left( - \frac{ 1}{r} \right) = \frac{( \cos\tilde\alpha_n - \cos\alpha_n) \pmb{\hat x}_n + (1 - \cos(\tilde\alpha_n-\alpha_n)) \pmb{\hat \xi}_n }{r_n^2\sin^2\alpha_n \xi_n} \ .
\ee
We substitute \ $\dis \pmb{\hat \xi}_n = \frac{\pmb{x}_{n+1} - \pmb{x}_n}{\xi_n r_n} $ from  \rf{defy1} and $\xi_n$  from  \rf{sinerule} and, after  elementary trigonometric calculations, we get 
\be  \label{rownosc}
 \bar \nabla^{\sss AVF} \left( - \frac{ 1}{r} \right) = \frac{( \pmb{x}_{n+1} \sin\tilde\alpha  + \pmb{x}_n \sin\alpha)  \sin\tilde\alpha  }{2 r_n^3 ( \cos \frac{\tilde\alpha_n-\alpha_n}{2})^2  \sin^2 \alpha} \ , 
\ee
and, using  \rf{differ} and \rf{sinerule}, we transform \rf{rownosc}  to the form \rf{Kep}.
\end{Proof}

\section{Locally exact modifications of  symmetric discrete gradient schemes }

Locally exact modification of the discrete gradient scheme
\be
   \pmb{y}_{n+1} = \pmb{y}_n + h_n  S {\bar \nabla}  H \ ,  \qquad    S = \mm 0 & I \\ - I & 0 \ema,
\ee
(where ${\bar \nabla}  H$ is a symmetric discrete gradient) is given by \cite{Ci-locex}:
\be  \label{lam}
  \pmb{y}_{n+1} = \pmb{y}_n + \Lambda_n {\bar \nabla}  H  \ , \qquad  \Lambda_n  = h_n  \tanhc \frac{h_n F'}{2} S \  ,
\ee
where  $F'$ is evaluated at $\pmb{\bar y}_n$ given by either $\pmb{\bar y}_n = \pmb{y}_n$  (LEX) or  $\pmb{\bar y}_n = \frac{1}{2} \left( \pmb{y}_n + \pmb{y}_{n+1} \right)$  (SLEX). Then,  $\tanhc (z) := z^{-1} \tanh (z)$ for $z\neq 0$ and $\tanhc (0) := 1$.  Tanhc is an analytic function:
\be
     \tanhc( x ) = 1 - \frac{1}{3} x^2  + \frac{2}{15} x^4 - \frac{17}{315} x^6 + \ldots 
\ee
Note that $\Lambda_n$ is of the same form for all  symmetric  discrete gradients, including the AVF discrete gradient.  
In the special case of separable Hamiltonians 
we have the following   useful theorem, see \cite{Ci-locex}. 

\ods
\begin{Th}   \label{Th-spec}
We assume $H (\pmb{x}, \pmb{p}) = T (\pmb{p}) + V (\pmb{x})$. Then  
the numerical scheme 
\be \ba{l}  \label{grad-del}
\pmb{\delta}_n^{-1} \left( {\pmb x}_{n+1} - {\pmb x}_n \right) = {\bar \nabla} T (\pmb{p}_n, \pmb{p}_{n+1})   ,  \\[2ex] \dis
(\pmb{\delta}_n^T)^{-1} \left( {\pmb p}_{n+1} - {\pmb p}_n \right) = - {\bar \nabla} V (\pmb{x}_n, \pmb{x}_{n+1})    , 
\ea \ee
(where $\bar \nabla$ is a symmetric discrete gradient), preserves exactly the energy integral (for any $h_n$-dependent  matrix \ $\pmb{\delta}_n$ such that $\pmb{\delta_n} \rightarrow h_n I$ for $h_n \rightarrow 0$) and is locally exact for
\be  \label{delta_n}
   \pmb{\delta}_n = h_n  \tanc \frac{h_n \Omega_n}{2} \ , \qquad \Omega_n^2  = T_{\pmb{p} \pmb{p}} (\pmb{\bar p}_n) V_{\pmb{x} \pmb{x}} (\pmb{\bar x}_n)  \ ,
\ee
where either $\pmb{\bar x}_n = \pmb{x}_n$, $\pmb{\bar p}_n = \pmb{p}_n$   or \  $\pmb{\bar x}_n = \frac{1}{2} \left( \pmb{x}_n + \pmb{x}_{n+1} \right)$, \ $\pmb{\bar p}_n = \frac{1}{2} \left( \pmb{p}_n + \pmb{p}_{n+1} \right)$.
\end{Th}
\ods

\begin{cor}
Locally exact  energy-preserving modifcations of the AVF method are  given by \rf{lam} for any $H$, and by \rf{grad-del} for  $H  = T (\pmb{p}) + V (\pmb{x})$.
\end{cor}

\ods

\section{Taylor expansions}  
\label{sec-Tay}

In this section we confine ourselves to the case \rf{H-p2}. 
Then, obviously,  $T_{\pmb{p} \pmb{p}}  = I$ and, as a consequence, $\Omega_n^2  =  V_{\pmb{x} \pmb{x}} (\pmb{\bar x}_n)$. Therefore
\be  \label{deltan}
\pmb{\delta}_n = h_n  \left(I + \frac{1}{12} h_n^2 V_{\pmb{x} \pmb{x}} + \frac{1}{120} h_n^4  V_{\pmb{x} \pmb{x}}^2  + 
\frac{17}{20 160} h_n^6  V_{\pmb{x} \pmb{x}}^3 +    \ldots \right) ,
\ee
where  $V_{\pmb{x x}}$  is evaluated at $\pmb{\bar x}_n$.

The system \rf{grad-del} can be (implicitly) solved with respect to  $\pmb{x}_{n+1}, \pmb{p}_{n+1}$. 
We are going to  derive $\pmb{x}_{n+1}, \pmb{p}_{n+1}$  in the form of  Taylor expansions in  $h$ (in this section we  denote $h_n \equiv h$):
\be  \label{ab}
 \pmb{x}_{n+1} = \pmb{x}_n + \sum_{k=1}^\infty \pmb{a}_k h^k \ , \qquad \pmb{p}_{n+1} = \pmb{p}_n + \sum_{k=1}^\infty \pmb{b}_k  h^k  \ ,
\ee 
where $\pmb{a}_k$, $\pmb{b}_k$ are vectors with components $a_{kj}$ and $b_{kj}$, respectively (they depend also on the index $n$, but we supress this dependence for the sake of clarity). 
Coefficients of the expansion \rf{deltan} will be denoted by  $c_k$ (matrices with entries $c_{k\mu\nu}$), i.e., 
\be  \label{delt}
  \pmb{\delta}_n  =  \sum_{k=1}^\infty h^k c_k ,
\ee
where, in particular, $c_1 = I$ (in order to assure the consistency of the modified scheme). We distinguish two cases. In the first case (LEX), we have
\be   \label{del-LEX}
  \pmb{\delta}_n^{\sss LEX} = h   I + \frac{1}{12} h^3 V_{\pmb{x}\pmb{x}} + \frac{1}{120} h^5 V_{\pmb{x}\pmb{x}}^2 + \ldots  ,
\ee
where $V_{\pmb{x x}}$ is evalueted at $\pmb{x}_n$. Note that  
$c_2 = c_4 = \ldots=0$, i.e.,  $\pmb{\delta}_n^{\sss LEX}$ is odd with respect to $h$. In the second case (SLEX) 
all derivatives on the right-hand side of \rf{deltan}  are  evaluated at $\frac{1}{2} (\pmb{x}_n + \pmb{x}_{n+1})$ and  $\pmb{x}_{n+1}$ is represented by the series \rf{ab}.  In particular,  
\be  \ba{l} \dis   \label{Vxx}
 V_{\pmb{x} \pmb{x}}  \left|_{\pmb{x} = \frac{1}{2} \left( \pmb{x}_n + \pmb{x}_{n+1} \right) } \right.   =  V_{\pmb{x} \pmb{x}}  +  V_{\pmb{x} \pmb{x} \pmb{x} } \left( \frac{1}{2} ( \pmb{x}_{n+1} + \pmb{x}_n)  - \pmb{x}_n \right) + \ldots  \\[1em] \dis
\qquad = V_{\pmb{x} \pmb{x}}   + \frac{1}{2} h V_{\pmb{x} \pmb{x} \pmb{x} } \pmb{a}_1 + \frac{1}{2} h^2  V_{\pmb{x} \pmb{x} \pmb{x} } \pmb{a}_2  + \frac{1}{8} h^2 V_{\pmb{x} \pmb{x} \pmb{x} \pmb{x} } (\pmb{a}_1, \pmb{a}_1) + \ldots 
\ea \ee
where all derivatives of $V$ on the right-hand side of \rf{Vxx} are evaluated at $\pmb{x}_n$. Thus
\be  \label{del-SLEX}
  \pmb{\delta}_n^{\sss SLEX}  =  h  I + \frac{1}{12} h^3 V_{\pmb{x}\pmb{x}} + \frac{1}{24} h^4   V_{\pmb{x} \pmb{x} \pmb{x}}  \pmb{a}_1  + \ldots  .
\ee
Note that $\pmb{\delta}_n$ is no longer odd in $h$.

Denoting \ $\mu$th component of $ \pmb{x}_{n+1} - \pmb{x}_n$ \ by \ $\Delta x^\mu$,  we expand \ $\gamma$th component of  \   $\bar \nabla V (\pmb{x}_n, \pmb{x}_{n+1})$    (shortly denoted by  $\bar \nabla_\gamma V$)  with respect to  $\Delta x^\mu$: 
\be \label{ABC}
\bar \nabla_\gamma V  = V,_\gamma  + \frac{1}{2} A_{\gamma\mu} \Delta x^\mu  +  \frac{1}{3!} B_{\gamma\mu\nu} \Delta x^\mu \Delta x^\nu +  \frac{1}{4!} C_{\gamma\mu\nu\sigma}   \Delta x^\mu   \Delta x^\nu  \Delta x^\sigma + \ldots 
\ee
 where  $A_{\gamma\mu}$, $B_{\gamma\mu\nu}$ and $C_{\gamma\mu\nu\sigma}$ are $\pmb{x}_n$-dependent coefficients (different for different discrete gradients) and summation over repeating Greek indices is assumed (from 1 to $m$). Here and in the following we denote:
\be
 V,_\gamma := \frac{\partial V}{\partial x^\gamma} \ , \quad  V,_{\gamma\mu} := \frac{\partial^2 V}{\partial x^\gamma \partial x^\mu} \ ,  \quad  \text{etc.}  
\ee
In the case of symmetric discrete gradients coefficients  $A_{\gamma\mu}$,  $B_{\gamma\mu\nu}$,  $C_{\gamma\mu\nu\sigma}$  have a special form.  

\begin{lem}
If the discrete gradient $\bar \nabla V$ is symmetric, then:
\be  \label{condit}
A_{\gamma\mu} = V,_{\gamma \mu} \ , \qquad  C_{\gamma \mu\nu\sigma} = 2 B_{\gamma \mu\nu},_\sigma - V,_{\gamma \mu\nu\sigma} .
\ee
\end{lem}

\begin{Proof} 
A symmetric discrete gradient satisfies \ $\bar \nabla_\gamma V (\pmb{x}_{n+1}, \pmb{x}_n)  =  \bar \nabla_\gamma V (\pmb{x}_n, \pmb{x}_{n+1}) $.  Taking into account that \ $\Delta x^\mu \rightarrow - \Delta x^\mu$ \  for \ $\pmb{x}_{n+1} \leftrightarrow \pmb{x}_n$, we compute
\be  \ba{l}  \dis  \label{V-odwrotnie}
\bar \nabla_\gamma V (\pmb{x}_{n+1}, \pmb{x}_n) = V,_\gamma (\pmb{x}_{n+1}) - \frac{1}{2} A_{\gamma\mu} (\pmb{x}_{n+1})  \Delta x^\mu     \\[1ex]\dis 
\qquad \qquad \quad  +  \frac{1}{3!} B_{\gamma\mu\nu} (\pmb{x}_{n+1})  \Delta x^\mu \Delta x^\nu -  \frac{1}{4!} C_{\gamma\mu\nu\sigma} (\pmb{x}_{n+1})  \Delta x^\mu   \Delta x^\nu  \Delta x^\sigma  + \ldots 
\ea \ee
Then, by standard Taylor expansions, we have
\[  \ba{l} \dis
V,_\gamma (\pmb{x}_{n+1}) = V,_\gamma + V,_{\gamma \mu} \Delta x^\mu + \frac{1}{2!} V,_{\gamma\mu\nu} \Delta x^\mu \Delta x^\nu + \frac{1}{3!} V,_{\gamma \mu\nu\sigma}  \Delta x^\mu   \Delta x^\nu   \Delta x^\sigma  + \ldots    \\[2ex]\dis
A_{\gamma\mu}  (\pmb{x}_{n+1}) = A_{\gamma\mu} + A_{\gamma\mu},_{\nu} \Delta x^\nu + \frac{1}{2!} A_{\gamma\mu},_{\nu\sigma} \Delta x^\nu   \Delta x^\sigma + \ldots    \\[2ex]\dis
B_{\gamma\mu\nu} (\pmb{x}_{n+1}) = B_{\gamma\mu\nu} + B_{\gamma\mu\nu},_{\sigma}  \Delta x^\sigma +  \frac{1}{2!} B_{\gamma\mu\nu},_{\sigma\rho} \Delta x^\sigma   \Delta x^\rho +  \ldots  \\[3ex]\dis
C_{\gamma\mu\nu\sigma} (\pmb{x}_{n+1}) = C_{\gamma\mu\nu\sigma} + C_{\gamma\mu\nu\sigma},_{\rho}  \Delta x^\rho +    \ldots  
\ea \]
where all quantities on the right-hand side are evaluated at $\pmb{x}_n$. Substituting these expansions into \rf{V-odwrotnie}, we get 
\be  \ba{l}  \dis
\bar \nabla_\gamma V (\pmb{x}_{n+1}, \pmb{x}_n) = V,_\gamma + \left( V,_{j\mu} - \frac{1}{2} A_{\gamma\mu} \right)   \Delta x^\mu     \\[2ex]\dis 
 \qquad \quad  +  \left(  \frac{1}{2!}  V,_{\gamma \mu\nu} -  \frac{1}{2!} A_{\gamma\mu},_\nu +  \frac{1}{3!} B_{\gamma\mu\nu} \right)  \Delta x^\mu \Delta x^\nu   \\[2ex]\dis   
 \qquad \quad  +   \left(  \frac{1}{3!}  V,_{\gamma\mu\nu\sigma} - \frac{1}{4} A_{\gamma\mu},_{\nu\sigma} + \frac{1}{3!} B_{\gamma\mu\nu},_\sigma  -  \frac{1}{4!} C_{\gamma\mu\nu\sigma}   \right)   \Delta x^\mu   \Delta x^\nu  \Delta x^\sigma  + \ldots 
\ea \ee
and comparing  these expressions with \rf{ABC}, term by term, we obtain \rf{condit}. 
\end{Proof}

\begin{lem} \label{lem-lex}
Coefficients $\pmb{a}_k$, $\pmb{b}_k$ for numerical solutions $\pmb{x}_{n+1}$, $\pmb{p}_{n+1}$ of \rf{grad-del}, where $\delta_n = \delta_n^{\sss LEX}$ is given by \rf{del-LEX}, have the following form (for  $k\leqslant 4$):
\be  \ba{l}   \label{Tay-grad}
a_{1\gamma} = p^\gamma \ ,   \qquad   b_{1\gamma} = - V,_\gamma \ , \\[1em]
a_{2\gamma} = - \frac{1}{2} V,_\gamma \  ,  \qquad   b_{2\gamma} = - \frac{1}{2} V,_{\gamma\mu} p^\mu   \ , \\[1em]
a_{3\gamma} = -  \frac{1}{6} V,_{\gamma\mu}  p^\mu  \   ,    \\[1em] 
b_{3\gamma} = \frac{1}{6} V,_{\gamma\mu}  V,_\mu -  \frac{1}{6} B_{\gamma\mu\nu} p^\mu  p^\nu       \   .
\ea \ee
\be \ba{l}  \label{Tay-grad-last}
  a_{4\gamma} =  \frac{1}{24} V,_{\gamma\mu}  V,_\mu -  \frac{1}{12} B_{\gamma\mu\nu} p^\mu  p^\nu   \ ,   \\[1em]
  b_{4\gamma} = \frac{1}{24}   V,_{\gamma\nu} V,_{\nu \mu}  p^\mu  + \frac{1}{12} \left(  B_{\gamma\mu\nu} + B_{\gamma\nu\mu} \right)  V,_\nu   p^\mu  \\[1ex] 
\qquad \qquad  \qquad \qquad   \qquad \qquad   - \frac{1}{24} \left( 2 B_{\gamma \mu \nu},_{\sigma } - V,_{\gamma\mu\nu\sigma} \right)  p^{\mu} p^{\nu} p^{\sigma}  .
\ea
\ee
\end{lem}

\begin{Proof}
Taking into account that \ $ \dis \Delta x^\mu =  \sum_{k=1}^\infty  a_{k\mu}  h^k$, we substitute  expansions \rf{ab}, \rf{delt} and \rf{ABC}  into \rf{grad-del}:
\be \ba{l}  \dis  \label{EQS}
   \sum_{k=1}^\infty  a_{k \gamma}  h^k   =  \left( \sum_{k=1}^\infty  c_{k \gamma \nu}  h^k  \right) \left( p^\nu     + \frac{1}{2}  \sum_{k=1}^\infty b_{k \nu}  h^k \right)  ,   \\[1em]\dis
 \sum_{k=1}^\infty  b_{k\gamma}  h^k = - \left(  \sum_{k=1}^\infty  c_{k\gamma \sigma}  h^k  \right) \left(    V,_\sigma + \frac{1}{2} A_{\sigma\mu} \Delta x^\mu  +  \frac{1}{3!} B_{\sigma\mu\nu} \Delta x^\mu \Delta x^\nu +  \ldots  \right)  ,
\ea \ee
where summation (from $1$ to $m$) over repeating Greek indices  is still  assumed.  
Equating coefficients by powers of $h$ we  compute (recurrently)  $a_{k\gamma}$ and $b_{k\gamma}$: 
\be  \ba{l}   \label{Tay-gr}
a_{1\gamma} = p^\gamma \ ,   \qquad   b_{1\gamma} = - V,_\gamma \ , \\[1em]
a_{2\gamma} = \frac{1}{2} b_{1\gamma} \  ,  \qquad  b_{2\gamma} = -  \frac{1}{2} A_{\gamma \mu} a_{1\mu}   \ , \\[1em]
a_{3\gamma} = \frac{1}{2} b_{2\gamma} + c_{3 \gamma\mu} p^\mu \   ,    \\[1em]
 b_{3\gamma} = -  \frac{1}{2} A_{\gamma\mu} a_{2\mu} -  \frac{1}{3!} B_{\gamma\mu\nu} a_{1\mu}  a_{1\nu} - c_{3\gamma\mu}  V,_\mu    \ ,   
\ea\ee
\be\ba{l} \label{Tay-gr-last}
a_{4\gamma} = \frac{1}{2} b_{3\gamma} + \frac{1}{2} c_{3 \gamma\mu} b_{1\mu} +  c_{4 \gamma\mu} p^\mu \ ,    \\[1em]
b_{4\gamma} = -  \frac{1}{2} A_{\gamma\mu} a_{3\mu} - \frac{1}{3!} B_{\gamma\mu\nu}  \left(  a_{1\mu} a_{2\nu} +   a_{2\mu}  a_{1\nu}  \right)   \\[1em] 
\qquad  - \frac{1}{4!} C_{\gamma\mu\nu\sigma} p^{\mu} p^{\nu} p^{\sigma}  -  \frac{1}{2} c_{3\gamma\nu} A_{\nu\mu} a_{1\mu} - c_{4 \gamma\nu} V,_\nu \ , \\[1ex]
\ea
\ee
where we assumed $c_{1\mu\nu} = \delta_{\mu\nu}$ (i.e., $c_1 = I$) and $c_{2\mu\nu}=0$, see  \rf{del-LEX}.  
 Then, we  use  \rf{condit} and take into account $c_{3 \gamma\mu} = \frac{1}{12} V,_{\gamma\mu}$ and  $c_{4\gamma\mu}=0$, see \rf{del-LEX}. Finally,  solving   equations \rf{Tay-gr} and \rf{Tay-gr-last} wherever possible, we obtain  \rf{Tay-grad} and \rf{Tay-grad-last}.  
\end{Proof}

\begin{lem} \label{lem-slex}
Coefficients $\pmb{a}_k$, $\pmb{b}_k$ for numerical solutions $\pmb{x}_{n+1}$, $\pmb{p}_{n+1}$ of \rf{grad-del}, where $\delta_n = \delta_n^{\sss SLEX}$ is given by \rf{del-SLEX}, have the form \rf{Tay-grad} for $k\leqslant 3$ and, for $k=4$ they are given by  
\be \ba{l}  \label{Tay-grad-last2}
  a_{4\gamma} =  \frac{1}{24} V,_{\gamma\mu}  V,_\mu -  \frac{1}{12} B_{\gamma\mu\nu} p^\mu  p^\nu   \ ,   \\[1em]
  b_{4\gamma} = \frac{1}{24}   V,_{\gamma\nu} V,_{\nu \mu}  p^\mu  + \frac{1}{12} \left(  B_{\gamma\mu\nu} + B_{\gamma\nu\mu} \right)  V,_\nu   p^\mu  \\[2ex] 
 \qquad   \qquad \qquad   - \frac{1}{24} V,_{\gamma\mu\nu} V,_{\nu} p^\mu - \frac{1}{24} \left( 2 B_{\gamma \mu \nu},_{\sigma } - V,_{\gamma\mu\nu\sigma} \right)  p^{\mu} p^{\nu} p^{\sigma}  .
\ea
\ee
\end{lem}

\begin{Proof}
 Similarly as in the proof of Lemma~\ref{lem-lex}, we substitute  expansions \rf{ab}, \rf{delt} and \rf{ABC}  into \rf{grad-del}, obtaining \rf{EQS},  \rf{Tay-gr} and \rf{Tay-gr-last}.  Then, we use  \rf{condit} and take into account  \rf{del-SLEX}, i.e.,  
\be \ba{l}
c_{3\gamma\mu} = \frac{1}{12} V,_{\gamma\mu} \ , \qquad   c_{4\gamma\nu} = \frac{1}{24} V,_{\gamma \nu\mu} a_{1\mu} = \frac{1}{24} V,_{\gamma \nu\mu} p^{\mu} \ . 
\ea \ee
As a result, we get \rf{Tay-grad} and \rf{Tay-grad-last2}. 
\end{Proof}

\begin{lem}
The Taylor expansion of the exact solution of the Hamiltonian system $\pmb{\dot p} = - V_{\pmb{x}}$,  $\pmb{\dot x} = \pmb{p}$ is given by
\be  \label{ab-exact}
 \pmb{x} (t+h)  = \pmb{x} (t)  + \sum_{k=1}^\infty \pmb{a}^{ex}_k h^k  , \qquad \pmb{p} (t+h)  = \pmb{p} (t) + \sum_{k=1}^\infty \pmb{b}^{ex}_k  h^k   ,
\ee 
where
\be  \ba{l}   \label{Tay-exact-ab}
a_{1\gamma}^{ex}  = p^\gamma \ ,   \qquad   b_{1\gamma}^{ex} = - V,_\gamma \ , \\[1em]
a_{2\gamma}^{ex} = - \frac{1}{2} V,_\gamma \  ,  \qquad   b_{2\gamma}^{ex} = - \frac{1}{2} V_{\gamma\mu} p^\mu   \ , \\[1em]
a_{3\gamma}^{ex} =  -  \frac{1}{6}  V,_{\gamma\mu} p^\mu  \   ,    \\[1em] 
b_{3\gamma}^{ex} =  \frac{1}{6}  V,_{\gamma\mu}  V,_\mu -  \frac{1}{6} V,_{\gamma\mu\nu} p^\mu  p^\nu       \   ,    \\[1em]
  a_{4\gamma}^{ex} = \frac{1}{24}  V,_{\gamma\mu}  V,_\mu -  \frac{1}{24} V,_{\gamma\mu\nu} p^\mu  p^\nu    \ ,   \\[1em]
  b_{4\gamma}^{ex} = \frac{1}{8} V,_{\gamma\mu\nu}  V,_\mu p^\nu + \frac{1}{24} V,_{\gamma\mu} V,_{\mu\nu} p^\nu    -  \frac{1}{24}  V,_{\gamma\mu\nu\sigma}   p^\mu p^\nu p^\sigma   \ .
\ea
\ee
\end{lem}

\begin{Proof}
Expanding $\pmb{y} (t+h)$, the exact solution to \rf{F(y)}, in a Taylor series, we get
\be  \ba{l} \label{exact-Tay}
\pmb{y} (t+h) = \pmb{y}(t) + h F + \frac{1}{2} h^2 F' F + \frac{1}{3!} h^3 \left( F''(F,F) + (F')^2 F \right)  \\[2ex] 
\quad \quad + \frac{1}{24}  h^4 \left(   F'''(F,F,F)  + 3 F''(F'F,F) + F' F''(F,F) + (F')^3 F   \right)  +  \ldots 
\ea  \ee  
where $F$ and its derivatives are evaluated at $\pmb{y} (t)$. 
In the case \rf{H-p2} we have
\be \ba{l}  \label{Tay-exact}
x^\gamma (t+h) = x^\gamma + h p^\gamma - \frac{1}{2} h^2 V,_\gamma - \frac{1}{3!} h^3  V,_{\gamma\mu} p^\mu 
  \\[1em]
\qquad \qquad + \frac{1}{4!} h^4 ( V,_{\gamma\mu} V,_\mu -  V,_{\gamma\mu\nu} p^\mu p^\nu )     + \ldots   \\[1em]
p^\gamma (t+h) =   p^\gamma - h V,_\gamma - \frac{1}{2!} h^2  V,_{\gamma\mu} p^\mu 
+ \frac{1}{3!} h^3 ( V,_{\gamma\mu} V,_\mu -  V,_{\gamma\mu\nu} p^\mu p^\nu )  \\[1em] 
\qquad \qquad + \frac{1}{4!} h^4 ( 3 V,_{\gamma\mu\nu} V,_\mu p^\nu + V,_{\gamma\mu} V,_{\mu\nu} p^\nu -  V,_{\gamma\mu\nu\sigma} p^\mu p^\nu p^\sigma ) + \ldots 
\ea \ee
which is equivalent to \rf{Tay-exact-ab}. 
\end{Proof}

\ods
\begin{Th} \label{Th-lex}
Locally exact modification (LEX) of a symmetric discrete gradient scheme defined by \rf{ABC} has the  order
\begin{itemize} 
\item at least 2 (always),
\item  at least 3 if and only if 
 $B_{\gamma\mu\nu} = V,_{\gamma\mu\nu}$,
\item at least 4 if and only if   $B_{\gamma\mu\nu} = V,_{\gamma\mu\nu} = 0$. 
\end{itemize}
\end{Th}

\begin{Proof}
It is sufficient to compare the expansion \rf{Tay-exact-ab} of the exact solution  with the expansion \rf{Tay-grad}, \rf{Tay-grad-last}.  We always  have: $a_{1\gamma} = a_{1\gamma}^{ex}$,   $b_{1\gamma} = b_{1\gamma}^{ex}$, $a_{2\gamma} = a_{2\gamma}^{ex}$,  $b_{2\gamma} = b_{2\gamma}^{ex}$,   $a_{3\gamma} = a_{3\gamma}^{ex}$. Then,  $b_{3\gamma} = b_{3\gamma}^{ex}$ only when $B_{\gamma\mu\nu} = V,_{\gamma\mu\nu}$. Finally, $a_{4\gamma} = a_{4\gamma}^{ex}$ and   $b_{4\gamma} = b_{4\gamma}^{ex}$ only if  $V,_{\gamma\mu\nu}=0$ (for any $\gamma, \mu, \nu$). 
\end{Proof}

\ods

\begin{Th}   \label{Th-slex}
Time reversible locally exact modification (SLEX) of a symmetric discrete gradient scheme \rf{ABC} has the order
\begin{itemize} 
\item at least 2 (always)  
\item at least  4 if and only if   $B_{\gamma\mu\nu} = V,_{\gamma\mu\nu}$,
\end{itemize}
\end{Th}

\begin{Proof} We compare the expansion \rf{Tay-exact-ab} of the exact solution  with the expansion \rf{Tay-grad}, \rf{Tay-grad-last2}.  Like in the proof of Theorem~\ref{Th-lex} we always have: $a_{1\gamma} = a_{1\gamma}^{ex}$,   $b_{1\gamma} = b_{1\gamma}^{ex}$, $a_{2\gamma} = a_{2\gamma}^{ex}$,  $b_{2\gamma} = b_{2\gamma}^{ex}$,   $a_{3\gamma} = a_{3\gamma}^{ex}$. Moreover,   $B_{\gamma\mu\nu} = V,_{\gamma\mu\nu}$  is a sufficient and necessary condition for  $b_{3\gamma} = b_{3\gamma}^{ex}$,  $a_{4\gamma} = a_{4\gamma}^{ex}$ and   $b_{4\gamma} = b_{4\gamma}^{ex}$.  
Note that for a symmetric discrete gradient  $B_{\gamma\mu\nu} = V,_{\gamma\mu\nu}$  implies $C_{\gamma\mu\nu\rho}=V,_{\gamma\mu\nu\rho}$, see  \rf{ABC}.
\end{Proof}

In order to compute coefficients $A_{\gamma\mu}, B_{\gamma\mu\nu}$ and  $C_{\gamma\mu\nu\rho}$  for the AVF method we expand the gradient $V,_\gamma$ (evaluated at  $\pmb{x}_n + t  \Delta \pmb{x}_n $)  in a Taylor series with respect to $\Delta \pmb{x}_n = \pmb{x}_{n+1} - \pmb{x}_n$. Thus we get
\be  \ba{l}  \dis
V,_\gamma (\pmb{x}_n + t  \Delta \pmb{x}_n )  = V,_\gamma + t V_{\gamma\mu} \Delta x^\mu  \\[1em] \dis \qquad 
 + \frac{1}{2} t^2 V_{\gamma\mu\nu} \Delta x^\mu \Delta x^\nu + \frac{1}{3!} t^3 V_{\gamma\mu\nu\sigma} \Delta x^\mu \Delta x^\nu \Delta x^\sigma + \ldots  
\ea \ee
Therefore, integrating term by term,
\be \ba{l}  \dis 
\bar \nabla^{\sss AVF}_\gamma  V = \int_0^1 V,_\gamma (\pmb{x}_n + t  \Delta \pmb{x}_n ) d t  \\[1em] \dis
\quad = V,_\gamma + \frac{1}{2}  V_{\gamma\mu} \Delta x^\mu + \frac{1}{3!} V_{\gamma\mu\nu} \Delta x^\mu \Delta x^\nu + \frac{1}{4!} V_{\gamma\mu\nu\sigma} \Delta x^\mu \Delta x^\nu \Delta x^\sigma + \ldots  
\ea \ee
which means that in the case of the AVF discrete gradient we have: 
\be
A_{\gamma\mu} = V,_{\gamma\mu} \ , \quad  B_{\gamma\mu\nu} = V,_{\gamma\mu\nu} \ , \quad  C_{\gamma\mu\nu\rho} = V,_{\gamma\mu\nu\rho} \ .
\ee

\begin{cor}
The locally exact method  AVF-LEX is of at least 3rd order  and the AVF-SLEX scheme  is of at least 4th order.  For special  potentials $V(\pmb{x})$ these orders can be higher, e.g., the AVF-LEX scheme is of 4th order for quadratic potentials. 
\end{cor}

\section{B-series}
\label{sec-B}

It is well known that the AVF method is a B-series method \cite{CMMOQW,QM}. The main reason is that the AVF discrete gradient has a Taylor expansion expressed in terms of the so-called elementary differentials, see \cite{HLW}. Indeed, denoting $\Delta \pmb{y}_{n} := \pmb{y}_{n+1} - \pmb{y}_n$,  expanding the right-hand side of \rf{AVF-y}  with respect to $\Delta \pmb{y}_{n}$, and integrating term by term, we get: 
\be  
\frac{\Delta \pmb{y}_{n}}{h} =   F + \frac{1}{2!} F' \Delta \pmb{y}_n + \frac{1}{3!} F'' ( \Delta \pmb{y}_n,  \Delta \pmb{y}_n) + \frac{1}{4!} F''' (\Delta \pmb{y}_n,  \Delta \pmb{y}_n, \Delta \pmb{y}_n) +  \ldots 
 \ee
(where $k$th derivative of $F$ is a vector-valued $k$-linear form).  
Treating this equation as an implicit equation for $\Delta \pmb{y}_n$ we can apply the standard fixed point method (starting from $\Delta \pmb{y}_n = 0$),  generating the formal series expression for the AVF scheme: 
\be
\Delta \pmb{y}_n = h F + \frac{1}{2} h ^2 F' F + \frac{1}{12} \left( 3 F' F' F + 2 F'' (F,F) \right) h^3 + \ldots 
\ee
We point out that at every iteration step we get a B-series (a sum of  elementary differentials). 

We will show that locally exact modifications preserve  this property:  both AVF-LEX  and AVF-SLEX  are B-series method. 
Locally exact AVF methods can be represented as
\be  \label{int}
\pmb{y}_{n+1} - \pmb{y}_n = \Lambda_n S^{-1} \int_0^1 F ((1-t) \pmb{y}_n + t 
\pmb{y}_{n+1} ) dt , 
\ee
where $\Lambda_n$ is defined by \rf{lam}. Note that
\be  \label{las}
\Lambda_n S^{-1} = h \left( I   - \frac{1}{12}  h^2 (F')^2 + \frac{1}{120} h^4 (F')^4 + \ldots  \right)  \ ,
\ee
where $F'$ is evaluated either at $\pmb{y}_n$ (AVF-LEX) or at $\frac{1}{2} (\pmb{y}_n + \pmb{y}_{n+1})$ (AVF-SLEX).

\begin{Th}
The AVF-LEX scheme \rf{int}, with $\Lambda_n S^{-1}$ evaluated at  $\pmb{y}_n$,  is a B-series method of the third order for any $F$. 
\end{Th}

\begin{Proof}
We  expand \rf{int} in the Taylor series with respect to $\Delta \pmb{y}_{n}$ and  integrate the right-hand side term by term, obtaining: 
\be  \ba{l} \dis
\Delta \pmb{y}_{n} = h \left( I   - \frac{1}{12}  h^2 (F')^2 + \frac{1}{120} h^4 (F')^4 + \ldots  \right) \\[1em] 
\dis \quad   \times \left(  F + \frac{1}{2!} F' \Delta \pmb{y}_n + \frac{1}{3!} F'' ( \Delta \pmb{y}_n,  \Delta \pmb{y}_n) + \frac{1}{4!} F''' (\Delta \pmb{y}_n,  \Delta \pmb{y}_n, \Delta \pmb{y}_n) +  \ldots \right)  
\ea \ee
where all derivatives are evaluated at $\pmb{y}_n$. 
The multiplication by powers of $F'$ transforms  elementary differentials into other elementary differentials. Therefore, applying the fixed point method in the identical way as in the AVF case, we express $\Delta \pmb{y}_n$ for the AVF-LEX scheme in terms of elementary differentials:  
\be \ba{l}  \label{AVF-LEX-Tay}
  \Delta \pmb{y}_n =  h F + \frac{1}{2} h ^2 F' F + \frac{1}{6} \left(  F' F' F +  F'' (F,F) \right) h^3  \\[2ex] 
 \quad \quad + \frac{1}{24} \left(   F'''(F,F,F)  + 4 F''(F'F,F) + 2 F' F''(F,F) + (F')^3 F   \right) h^4  +  \ldots 
\ea  \ee  
Finally, comparing  \rf{AVF-LEX-Tay} with \rf{exact-Tay}, we see that AVF-LEX is, at least, of the third order. 
\end{Proof}

\begin{Th}
The AVF-SLEX scheme, given by \rf{int} with $\Lambda_n S^{-1}$ evaluated at $\frac{1}{2} (\pmb{y}_n + \pmb{y}_{n+1})$,   is a B-series method of the 4th order (at least) for any $F$. 
\end{Th}

\begin{Proof}
We have to evaluate $\Lambda_n S^{-1}$, given by \rf{las}, at  $\frac{1}{2} (\pmb{y}_n + \pmb{y}_{n+1})$, i.e., at $\pmb{y}_n + \frac{1}{2} \Delta \pmb{y}_n$. We denote  $ \widehat{F'} :=   F'(\pmb{y}_n + \frac{1}{2} \Delta \pmb{y}_n)$. From \rf{int} we obtain:
\be  \ba{l} \dis  \label{rownosc2}
\Delta \pmb{y}_{n} = h \left( I   - \frac{1}{12}  h^2 (\widehat{F'})^2 + \frac{1}{120} h^4 (\widehat{F'})^4 + \ldots  \right) \\[1em] 
\dis \quad   \times \left(  F + \frac{1}{2!} F' \Delta \pmb{y}_n + \frac{1}{3!} F'' ( \Delta \pmb{y}_n,  \Delta \pmb{y}_n) + \frac{1}{4!} F''' (\Delta \pmb{y}_n,  \Delta \pmb{y}_n, \Delta \pmb{y}_n) +  \ldots \right)  .
\ea \ee
Expanding $\widehat{F'}$ with respect to $\Delta \pmb{y}_n$, we get
\be 
 \widehat{F'} \equiv  F'(\pmb{y}_n + \frac{1}{2} \Delta \pmb{y}_n) = F' + \frac{1}{2} F''(\Delta \pmb{y}_n, \cdot) + \frac{1}{8} F'''(\Delta \pmb{y}_n, \Delta \pmb{y}_n, \cdot) + \ldots  
 \ee
Suppose that $\Delta \pmb{y}_n$ and $\pmb{w}$ are given by some B-series. Then 
\be
 \widehat{F'} \ \pmb{w} = F' \pmb{w}  + \frac{1}{2} F''(\Delta \pmb{y}_n, \pmb{w}) + \frac{1}{8} F'''(\Delta \pmb{y}_n, \Delta \pmb{y}_n, \pmb{w}) + \ldots
\ee 
is a B-series as well, see \cite{HLW}. Therefore,  $(\widehat{F'})^n \pmb{w}$ is also a B-series for any $n \in \N$. Applying the fixed point method for solving \rf{rownosc2}, we obtain B-series at every iteration step.  Thus we can  express $\Delta \pmb{y}_n$ for the AVF-SLEX scheme in terms of elementary differentials:  
\be  \ba{l} \dis \label{AVF-SLEX-Tay}
  \Delta \pmb{y}_n =  h F + \frac{1}{2} h ^2 F' F + \frac{1}{6} \left(  F' F' F +  F'' (F,F) \right) h^3  \\[2ex]  
 \quad \quad + \frac{1}{24} \left(   F'''(F,F,F)  + 3 F''(F'F,F) + F' F''(F,F) + (F')^3 F   \right) h^4  +  \ldots 
\ea  \ee  
All these terms (up to 4th order) are identical with the  expansion of the exact solution. Therefore, AVF-SLEX is of the 4th order. 
\end{Proof}

\section{Numerical experiments}
\label{sec-num}

We tested numerically locally exact AVF schemes  on spherically symmetric potentials (circular orbits for  the Coulomb-Kepler problem  and an anharmonic oscillator). Our algorithms use  simplified expressions for the AVF  discrete gradients derived in Section~\ref{sec-AVF}.  The function $\tanhc$ of  a matrix is computed in the standard way  (by diagonalization of the matrix) as any other analytic function. 
This straightforward method is expensive for large $m$ and  then one may use  special techniques for computing exponential integrators \cite{HL,NW}.

We took into account that locally exact modifications are more expensive. Therefore, in numerical experiments we use different time steps in order to get  the same computational costs. The cost was estimated by the number of function evaluations. 

\begin{figure}
\caption{Error of numerical solutions (at $t = 100$) for a circular orbit ($R=3.5$, exact period $T\approx 41.14$) in the Coulomb-Kepler potential  $V (r) = -1/r$ as a function of the scaled time step $\tilde h$.  AVF: white discs ($\tilde h = h$),  AVF-LEX:  black discs ($\tilde h/h \approx 1.8, 2.0, 1.7, 1.9$),  AVF-SLEX: black squares ($\tilde h/h = 2.8, 3.3, 3.2, 3.1$). } 
 \label{Kepler}  \par   \ods
\includegraphics[width=\textwidth]{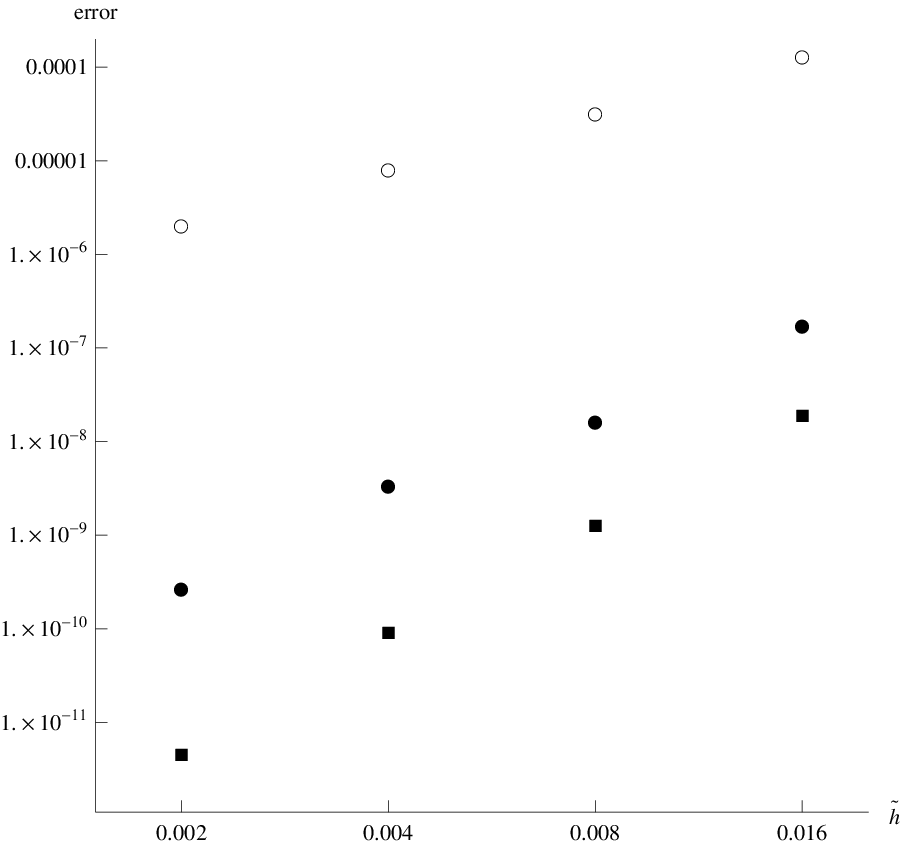} \par
\end{figure}

\begin{figure}
\caption{Error of numerical solutions (at $t =100$) for a circular orbit ($R=3.5$, exact period $T\approx 8.80$) in the potential of an anharmonic oscillator $V (r) = 0.5 r^2  - 0.01 r^4$ as a function of the scaled time step $\tilde h$.   AVF: white discs  ($\tilde h = h$),  AVF-LEX:  black discs ($\tilde h/h \approx 1.9, 1.6, 1.7, 1.7$),  AVF-SLEX: black squares ($\tilde h/h \approx 2.9, 2.8, 2.8, 2.9$). } 
 \label{oscyl}  \par   \ods
\includegraphics[width=\textwidth]{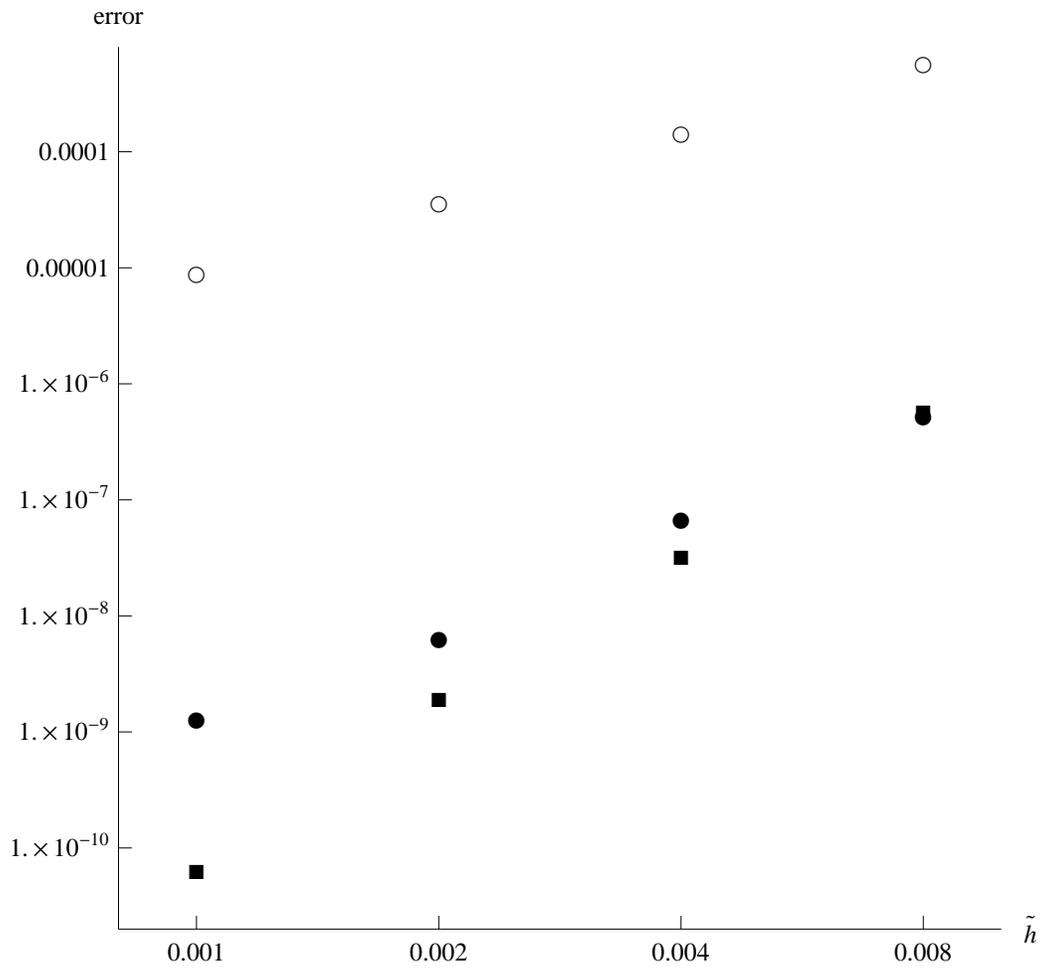} \par
\end{figure}

Results of numerical experiments are very promising, see Fig.~\ref{Kepler} and Fig.~\ref{oscyl}. 
The scaled time step $\tilde h $  equals $h$ for the AVF scheme. In the case of  AVF-LEX the actual time step $h$ is  about  2 times greater than $\tilde h$ (for more precise values see  figure captions), and for AVF-SLEX this factor is about 3.   In spite of this handicap locally exact modifications are more accurate by several orders of magnitude in comparison with the standard AVF scheme.

Note that the $h$-dependence of the error is practically linear (in both cases). The slopes of these lines are given by $2.00, 3.02, 3.99$ (Fig.~\ref{Kepler}) and $2.00, 2.95, 4.36$ (Fig.~\ref{oscyl}) which is in a good agreement with computed orders of the corresponding schemes. Note that for larger values of $\tilde h$  the accuracy of AVF-LEX can be greater than the accuracy of AVF-SLEX.

The higher order of  AVF-LEX and AVF-SLEX (as compared to AVF) explains only in part  their high accuracy. We point out that locally exact schemes are very accurate in a neighbourhood (not very small, in fact) of the stable equilibrium \cite{Ci-CMA,Ci-locex-PLA}.  Local exactness  has also some additional advantages \cite{Ci-CMA} and is worthwhile to be studied more carefully.

\ods

\no {\bf Acknowledgments}. I am grateful to anonymous reviewers for corrections and suggestions which improved the paper. Research supported in part  by the National Science Centre (NCN) grant no. 2011/01/B/ST1/05137.

\end{document}